\documentclass[letterpaper,10pt,conference]{ieeeconf}
\IEEEoverridecommandlockouts
\usepackage{graphicx}
\usepackage{subcaption}
\usepackage{amsmath}
\usepackage{amsfonts}
\usepackage{amssymb}
\usepackage{bm}
\usepackage{mathtools}
\usepackage{hyperref}
\usepackage{xcolor}
\usepackage{booktabs}

\newtheorem{assumption}{Assumption}
\newtheorem{theorem}{Theorem}
\newtheorem{remark}{Remark}

\newcommand{\X}{\mathcal X}
\newcommand{\U}{\mathcal U}
\newcommand{\Kcal}{\mathcal K}
\newcommand{\E}{\mathbb E}
\newcommand{\R}{\mathbb R}

\newcommand{\argmin}{\operatorname*{arg\,min}}

\title{\LARGE \bf Dual Model Predictive Control \\ through Sampling-Based Belief Propagation}

\author{Tren Baltussen$^{1}$, Johannes K\"ohler$^{2}$, Alexander Katriniok$^{1}$, Maurice Heemels$^{1}$%
\thanks{$^{1}$T. Baltussen, A. Katriniok and M. Heemels are with the Department of
Mechanical Engineering, Eindhoven University of Technology, the Netherlands.
$^{2}$J. K\"ohler is with the Department of Mechanical Engineering, Imperial College
London, UK. Email: {\tt t.m.j.t.baltussen@tue.nl}.}}

\begin{document}

\maketitle
\thispagestyle{empty}
\pagestyle{empty}

\begin{abstract}
    We propose a sampling-based dual model predictive control (MPC) approach that combines belief propagation with parameterized feedback policy optimization. By sampling uncertain dynamics and process noise, we approximate the evolution of the hyperstate over the prediction horizon, allowing future control decisions to adapt to anticipated observations and belief updates. This enables the controller to account for the future value of information without explicit exploration incentives. We establish almost-sure uniform convergence of the sample-average objective to the finite-horizon expected objective over a compact policy class. The approach is specialized to Bayesian linear regression and evaluated on stochastic systems with constant and time-varying uncertain dynamics. We compare our approach to certainty equivalent control, and an MPC that optimizes over open-loop control sequences rather than feedback policies, demonstrating the benefits of our method in terms of both control performance and information gain.

\end{abstract}

\section{Introduction}

Dual control addresses the problem of simultaneously regulating a dynamical system and learning its uncertain dynamics. As originally recognized by Feldbaum \cite{feldbaum1960dual}, control actions influence both the physical state and the information available for future decisions. The optimal Bayesian control policy accounts for this dual effect by balancing immediate control performance against the future benefits of improved model knowledge \cite{barshalom1974dual,astrom1989adaptive}. However, computing this policy requires solving a stochastic dynamic program over the posterior distribution of the physical state and of the unknown dynamics, the so-called \emph{hyperstate}, which is generally computationally intractable \cite{mesbah2018dual}.

Model predictive control (MPC) can provide a tractable approximation of the optimal control problem, typically by an open-loop input sequence over a finite prediction horizon \cite{rawlings2017mpc}. 
By optimizing over feedback policies and closed-loop predictions the controller can anticipate how current decisions affect future observations, belief updates, and subsequent control decisions. This enables the value of information to influence the control policy, thereby capturing the dual effect. Motivated by this idea, we propose a tractable sampling-based dual MPC approach that optimizes parameterized feedback policies over sampled state and belief trajectories.

\subsection{Related Work}
Approximate dual control methods are commonly classified as explicit or implicit \cite{mesbah2018dual}.
Explicit approaches encourage learning through exploration terms, excitation signals, or learning constraints \cite{soloperto2020augmenting,soloperto2023guaranteed,baltussen2026dual}, but typically require additional design choices beyond the original control objective.

Implicit approaches instead approximate the underlying stochastic optimal control problem to account for the future value of information. These dual-control approximations include, e.g., predicted estimation errors and analytic uncertainty propagation \cite{heirung2017dual,tse1973}, or approximate Bayesian control \cite{klenske2016dual}.
Various sampling-based MPC approaches have been proposed to approximate the dual control problem. Here, the multi-stage decision problem is formulated as a scenario tree, where different `branches' of the tree use individual control inputs to account for different uncertainty realizations, allowing future decisions to adapt to anticipated observations \cite{arcari2020adp,arcari2020dual, THANGAVEL2018}
However, the number of decision variables grows exponentially with the number of scenarios and the prediction horizon, which may pose challenges to the practical applicability of these methods.

\subsection{Contributions}
We propose a sampling-based dual MPC formulation that optimizes a shared sequence of parameterized hyperstate-feedback policies.
In contrast with other sampling-based approaches, the computed policy is actually an implementable feedback law that can be rolled out for any uncertainty realization and does not rely on hypothetical realizations.
By propagating sampled state and posterior trajectories, the formulation accounts for anticipated observations and their influence on future control decisions, thereby intrinsically enabling the dual effect.
Moreover, the number of policy parameters is independent of the number of samples.
We establish almost-sure uniform convergence of the sample-average objective to the exact finite-horizon Bayesian objective within the chosen policy class under standard regularity assumptions.
The approach is specialized to Bayesian linear regression and evaluated on stochastic examples with time-invariant and time-varying uncertain input gains.

\section{Finite-Horizon Optimal Control Problem}
\label{sec:exact_problem}
Consider the discrete-time system
\begin{equation}
    x_{k+1}=g(x_k,u_k)+w_k,
    \label{eq:general_dynamics}
\end{equation}
where $x_k\in\X=\R^{n_x}$ is the fully observed state, $u_k\in\U\subseteq\R^{n_u}$ is the control input at time index $k \in \mathbb{N}_0$, and $g \in \mathcal{G}: \X\times\U\to\R^{n_x}$ is a continuous unknown, time-invariant transition function.
The dynamics $g$ is a random function with a known prior distribution $\mathbb{P}_g$.
The disturbances $w_k\sim\mathbb P_w$ are independent and identically distributed (i.i.d.), with known distribution, and independent of $g$. Time-varying dynamics are discussed in
Remark~\ref{rem:time_varying_dynamics}.

Let $\gamma_k=\mathbb P(g\mid\mathcal I_k)$ denote
the posterior distribution of $g$, given the information available at time $k$,
$\mathcal I_k=(x_0,u_0,\ldots,u_{k-1},x_k)$.
We define the current hyperstate as the joint distribution of the physical state and the posterior distribution of the dynamics:
\(
    \xi_k=\delta_{x_k}\otimes\gamma_k \in\Xi,
\)
where $\Xi:=\mathcal P(\X\times\mathcal G)$ denotes the space of probability measures on $\X\times\mathcal G$.
Since the state is fully observed at time $k$,
we identify the hyperstate with the pair
\begin{equation}
    \xi_k=(x_k,\gamma_k).
\end{equation}
while the future hyperstate is random due to the uncertain dynamics and disturbances.
Using Bayes rule and independence of the disturbances, the posterior distribution $\gamma_{k+1}=\mathbb{P}( g \mid \mathcal{I}_{k+1})$ can be expressed as
\begin{equation}
    \gamma_{k+1}
    =\mathcal B(\gamma_k,x_k,u_k,x_{k+1}),
    \label{eq:bayes_update}
\end{equation}
i.e., depending only on $\gamma_k$, on the information at time $k$, and the new information from one transition $(x_k,u_k,x_{k+1})$.

Given an initial hyperstate $\xi$, we consider the $N$-step finite-horizon stochastic optimal control problem
\begin{equation}
\begin{aligned}
    \inf_{\boldsymbol\pi\in\Pi^N}\quad
    &\E_{\boldsymbol\omega}
    \left[
        \sum_{i=0}^{N-1}\ell(x_i,u_i)
        +\ell_N(x_N)
        \,\middle|\,\xi_0=\xi
    \right]\\
    \mathrm{s.t.}\quad
    &u_i=\pi_i(\xi_i),\\
    &x_{i+1}=g(x_i,u_i)+w_i,\\
    &\gamma_{i+1}
      =\mathcal B(\gamma_i,x_i,u_i,x_{i+1}),\\
    &\xi_i=(x_i,\gamma_i),\\
    &\xi_0=\xi,\qquad i=0,\ldots,N-1,
\end{aligned}
\label{eq:exact_bayesian_dual_ocp}
\end{equation}
where $\boldsymbol\pi=(\pi_0,\ldots,\pi_{N-1}) \in \Pi^N$
is a sequence of causal feedback laws
$\pi_i \in \Pi :\Xi\to\U$,
and $\ell:\X\times\U\to\R_{\geq0}$ and
$\ell_N:\X\to\R_{\geq0}$ are the stage and
terminal costs.
Here, the expectation over the uncertainty $\boldsymbol\omega=(g,w_0,\ldots,w_{N-1})$ is conditioned on the initial hyperstate $\xi$.
Because control actions influence both future states and observations, they also affect subsequent posterior distributions and feedback decisions. The optimization therefore intrinsically accounts for the value of information, giving rise to the dual effect.

\section{Sampling-based Belief Propagation}
\label{sec:sampling_based_belief_propagation}
As the exact problem \eqref{eq:exact_bayesian_dual_ocp} is generally intractable \cite{mesbah2018dual}, we next restrict the policy class $\Pi$ and approximate the expectation over $\boldsymbol\omega$ by a sample-average approximation (SAA) \cite{shapiro2014lectures}.

\subsection{Policy Class}
We restrict the admissible feedback laws to a
finite-dimensional parametrized class of the form
\begin{equation}
    \Pi_{\Kcal}
    :=\{\pi_K:K\in\Kcal\},
    \qquad
    \pi_K:\Xi\to\U,
\end{equation}
where the policy parameter set $\Kcal\subset\R^n$ is compact. Then, we optimize the policy parameter sequence $\mathbf K=(K_0,\ldots,K_{N-1})\in\Kcal^N$ over a prediction horizon $N$, with the feedback policy $u_k=\pi_{K_k}(\xi_k)$.

For a given uncertainty realization $\boldsymbol\omega$ over the prediction horizon, an initial hyperstate $\xi$, and a policy sequence $\mathbf{K}$, the corresponding closed-loop rollout cost is
\begin{equation}
    C_N(\mathbf K,\boldsymbol\omega;\xi)
    :=
    \sum_{k=0}^{N-1}
    \ell\left(x_k,\pi_{K_k}(\xi_k)\right)
    +\ell_N(x_N),
    \label{eq:rollout_cost}
\end{equation}
where the hyperstate $(\xi_k)_{k=0}^{N-1}$ is propagated
according to the dynamics \eqref{eq:general_dynamics} and Bayesian
update \eqref{eq:bayes_update}, starting from $\xi_0$.

The resulting finite-horizon objective is
\begin{equation}
    J_N(\mathbf K;\xi)
    :=
    \E_{\boldsymbol\omega}
    \left[C_N(\mathbf K,\boldsymbol\omega;\xi)\right].
    \label{eq:policy_class_objective}
\end{equation}

\begin{assumption}[Regularity]
\label{ass:saa_regular}
For any fixed hyperstate $\xi\in\Xi$, suppose that:
\begin{enumerate}
    \renewcommand{\labelenumi}{\roman{enumi})}
    \item The parameter set $\Kcal^N$ is nonempty
    and compact.
    \item The rollout cost
    $C_N(\mathbf K,\boldsymbol\omega;\xi)$ is continuous
    in $\mathbf K$ for almost every $\boldsymbol\omega$.
    \item The finite-horizon rollout costs satisfy
    \begin{equation}
        \E_{\boldsymbol\omega}\left[
        \sup_{\mathbf K\in\Kcal^N}
        |C_N(\mathbf K,\boldsymbol\omega;\xi)|
        \right]<\infty.
    \end{equation}
    \item The uncertainty realizations
    $\boldsymbol\omega^{(1)}, \boldsymbol\omega^{(2)}\ldots,
    $ are independent and
    identically distributed according to their
    conditional distribution given $\xi$.
\end{enumerate}
\end{assumption}

Condition (i) holds for a finite-dimensional policy
parameterization with closed and bounded parameter sets.
Condition (ii) holds assuming continuous dynamics, Bayesian update, feedback policies, and cost functions. Such compactness and continuity requirements are standard in stochastic MPC \cite[Sec.~II-C]{mcallister2023nonlinear}.
Condition (iii) holds, e.g., for the examples considered in Section~\ref{sec:numerics}, due to the dynamics, finite-prediction horizon, bounded policy coefficients, and quadratic costs. Finally, condition (iv) follows directly from independently sampling the complete uncertainty trajectories.

\subsection{Sample-Average Approximation}
\label{sec:saa}
We approximate the expectation in \eqref{eq:policy_class_objective} using samples of the uncertainty $\boldsymbol{\omega}$, yielding a SAA \cite{shapiro2014lectures} that can be optimized over the policy parameter sequence $\mathbf K$. By propagating the sampled hyperstate under parameterized feedback policies, the optimization accounts for anticipated belief updates and their influence on future control decisions, thereby capturing the dual effect.

Let
\(
    \boldsymbol\omega^{(1)},\ldots,\boldsymbol\omega^{(S)}
\)
be $S$ independent samples from the conditional distribution of
\(\boldsymbol\omega\) given \(\xi\). The SAA of \eqref{eq:exact_bayesian_dual_ocp} is
\begin{equation}
    \widehat{\mathbf K}_{S}
    \in
    \argmin_{\mathbf K\in\Kcal^N}
    \widehat J_{N,S}(\mathbf K;\xi),
    \label{eq:saa_problem}
\end{equation}
where
\begin{equation}
    \widehat J_{N,S}(\mathbf K;\xi)
    :=
    \frac1S
    \sum_{s=1}^S
    C_N(\mathbf K,\boldsymbol\omega^{(s)};\xi).
    \label{eq:saa_objective}
\end{equation}
Under Assumption~\ref{ass:saa_regular}, the minimum of the sampled problem \eqref{eq:saa_problem} is attained almost surely on $\Kcal^N$.
The same samples are used for all candidate policy sequences during the optimization.
Thus, \eqref{eq:saa_problem} is a deterministic optimization problem that can be solved with standard optimization methods.
The policy parameters $\mathbf K$ are shared across all scenarios, whereas the predicted hyperstates evolve according to their individual sampled observations.
Thus, while the sampled parameter remains fixed, the predicted posterior and subsequent control actions depend on the anticipated observations of the sampled hyperstate.

\subsection{Finite-Horizon SAA Guarantees}
\label{sec:theory}
Under mild regularity conditions, we establish almost-sure uniform convergence of the sampled objective to the exact expected objective over the compact policy-parameter set using standard results from SAA theory \cite{shapiro2014lectures}.

\begin{theorem}[Finite-horizon SAA consistency]
\label{thm:saa_consistency}
Fix a hyperstate \(\xi\), and suppose that
Assumption~\ref{ass:saa_regular} holds. Then
\begin{equation}
    \sup_{\mathbf K\in\Kcal^N}
    \left|
        \widehat J_{N,S}(\mathbf K;\xi)
        -
        J_N(\mathbf K;\xi)
    \right|
    \to 0
    \quad
    \text{almost surely}
    \label{eq:uniform_saa_convergence}
\end{equation}
as \(S\to\infty\), i.e., the sampled objective converges uniformly over the policy-parameter set to the exact expected objective.
\end{theorem}

\begin{proof}
By Assumption~\ref{ass:saa_regular}, the sample-average functions satisfy the
conditions of the uniform law of large numbers on the compact set \(\Kcal^N\):
the integrand is almost surely continuous in \(\mathbf K\), dominated by an
integrable random variable, and the samples are i.i.d.. Hence
\eqref{eq:uniform_saa_convergence} follows from the uniform law of large numbers
\cite[Sec.~7.2.5, Thm.~7.48]{shapiro2014lectures}.
\end{proof}

The theorem establishes uniform convergence of the sampled objective to the exact expected objective over the chosen policy class. Consequently, the expected cost attained by any global minimizer of the sampled problem converges to the minimum expected cost within the chosen policy class $\Kcal^N$ as \(S\to\infty\). This does not imply convergence to the unrestricted Bayesian optimal controller to \eqref{eq:exact_bayesian_dual_ocp}. 
Under Assumption~\ref{ass:saa_regular}, the sampling error converges to zero almost surely, while the remaining difference to \eqref{eq:exact_bayesian_dual_ocp} depends on the expressiveness of the chosen feedback policies.

\begin{remark}[Time-varying dynamics]
\label{rem:time_varying_dynamics}
The framework can be extended to time-varying dynamics of the form
\begin{equation}
\begin{aligned}
    x_{k+1} &= g(x_k,u_k;\theta_k)+w_k,\\
    \theta_{k+1} &= h(\theta_k,v_k),
\end{aligned}
\end{equation}
where $\theta_k\in\Theta$ is an unknown parameter, $h:\Theta\times\mathcal V\to\Theta$ is a known measurable function, and $v_k\in\mathcal V$ is an i.i.d. disturbance with known distribution $\mathbb P_v$, independent of $w_k$.

With a slight abuse of notation, the belief $\gamma_k$ now describes the
distribution of $\theta_k$. The Bayesian update consists of two steps. First, the current belief is conditioned on the new observation
$x_{k+1}$
\begin{equation}
\label{eq:time_varying_bayes_update}
\gamma_k^+(d\theta)
\propto
p_w(x_{k+1}-g(x_k,u_k;\theta))
\,\gamma_k(d\theta),
\end{equation}
where $p_w$ is the known density of $w_k$.

Second, the conditioned belief is predicted
through the parameter dynamics:
\begin{equation}
\label{eq:time_varying_bayes_prediction}
\gamma_{k+1}(A)
=
\int_\Theta
\mathbb P(h(\theta,v_k)\in A)
\,\gamma_k^+(d\theta),
\end{equation}
for any measurable $A\subseteq\Theta$.
Together, the conditioning \eqref{eq:time_varying_bayes_update} and prediction \eqref{eq:time_varying_bayes_prediction} steps define the Bayesian update map $\mathcal B$. 
Theorem~\ref{thm:saa_consistency} also applies
to the time-varying case as the scenarios remain independent, while the parameters evolve according to $h$ within each scenario.
The time-invariant case is recovered when $h(\theta,v)=\theta$.
\end{remark}

The approach requires a Bayesian update map $\mathcal B$ and the ability to sample model realizations from the current posterior $\gamma_k$.
We expect that similar constructions can be derived for more general stochastic processes on measurable function spaces, such as Gaussian-process models \cite{rasmussen2006gaussian}. These extensions are left for future work.

\section{Special Case: Bayesian Linear Regression}
\label{sec:finite_models}
In this section, we show how the proposed approach simplifies in the special case of Bayesian linear regression.

\subsection{Bayesian Linear Regression}
\label{sec:bayesian_linear_regression}
Consider the parametric transition model
\begin{equation}
    x_{k+1}=g(x_k,u_k;\theta)+w_k,\qquad
    g(x,u;\theta)=\Phi(x,u)\theta,
    \label{eq:blr_transition_model}
\end{equation}
where $w_k\sim\mathcal N(0,W)$ is Gaussian with covariance $W\succ0$ and the unknown parameter $\theta\in\mathbb R^{n_\theta}$ is fixed over time. We consider a Gaussian prior for $\theta \sim\mathcal N(\mu_0,\Sigma_0)$.
Given the information
$\mathcal I_k$, we obtain
\(
    \theta\mid\mathcal I_k\sim\mathcal N(\mu_k,\Sigma_k),
    \Sigma_k\succ0.
\)
Introducing the information variables
\(
    \Lambda_k:=\Sigma_k^{-1},    \eta_k:=\Lambda_k\mu_k,
\)
gives the hyperstate $\xi_k=(x_k,\eta_k,\Lambda_k)$. Given a transition
observation $(x_{k},u_k,x_{k+1})$, the Bayesian update is
\begin{equation}
\begin{aligned}
    \Lambda_{k+1}&=\Lambda_k+\Phi(x_k,u_k)^\top W^{-1}\Phi(x_k,u_k),\\
    \eta_{k+1}&=\eta_k+\Phi(x_k,u_k)^\top W^{-1}x_{k+1}.
\end{aligned}
\label{eq:blr_update_map}
\end{equation}
The corresponding posterior moments are
\(
    \Sigma_{k+1}=(\Lambda_{k+1})^{-1},
\)
\(
    \mu_{k+1}=\Sigma_{k+1}\eta_{k+1}.
\)
To sample the unknown parameter, let
$z^{(s)}\sim\mathcal N(0,I_{n_\theta})$ be an independent standard Gaussian random variable of dimension $n_\theta$. With $L_kL_k^\top=\Sigma_k$, the model realizations are obtained as
$\theta_k^{(s)} = \mu_k+L_kz^{(s)}$ for $s=1,\ldots,S$.
The samples $\left(z^{(s)}\right)$ remain fixed throughout the prediction horizon. Starting from $\xi_{0}^{(s)}=\xi$, the sampled dynamics are
\begin{equation}
    x_{k+1}^{(s)}
    =
    \Phi(x_{k}^{(s)},u_{k}^{(s)})
    \theta_k^{(s)}
    +w_{k}^{(s)},
    \label{eq:blr_sampled_state_recursion}
\end{equation}
where
$u_{k}^{(s)}
=\pi_{K_k}(\xi_{k}^{(s)})$.
After each predicted transition, the posterior is
updated according to \eqref{eq:blr_update_map},
giving the next hyperstate
$\xi_{k+1}^{(s)}$.

The resulting sample-average MPC problem for Bayesian linear regression at a given hyperstate $\xi$ is
\begin{equation}
\begin{aligned}
\widehat{\mathbf K}_{S}(\xi)
\in \argmin_{\mathbf K\in\Kcal^N}
\quad &
\frac{1}{S}\sum_{s=1}^S
\left[
\sum_{k=0}^{N-1}
\ell(x_{k}^{(s)},u_{k}^{(s)})
+\ell_N(x_{N}^{(s)})
\right] \\
\mathrm{s.t.}\quad
&u_{k}^{(s)}
=\pi_{K_k}(\xi_{k}^{(s)}),\\
&x_{k+1}^{(s)}
=\Phi(x_{k}^{(s)},u_{k}^{(s)})
\theta_k^{(s)}+w_{k}^{(s)},\\
&\eta_{k+1}^{(s)}
=\eta_k^{(s)}
+\Phi_k^{(s)\top}W^{-1}x_{k+1}^{(s)},\\
&\Lambda_{k+1}^{(s)}
=\Lambda_k^{(s)}
+\Phi_k^{(s)\top}W^{-1}\Phi_k^{(s)},\\
&(x_{0}^{(s)}, \eta_{0}^{(s)},\Lambda_{0}^{(s)})
=\xi,\\
&k=0,\ldots,N-1,\quad s=1,\ldots,S.
\end{aligned}
\label{eq:blr_saa_ocp}
\end{equation}

\section{Receding-Horizon Control}
\label{sec:rhc_implementation}
To apply the finite-horizon formulation to a control task with duration larger than $N$, the sampled problem is solved in a receding-horizon fashion \cite{rawlings2017mpc}.
At time $k$, the measured state and current model posterior form the hyperstate $\xi_k=(x_k,\gamma_k)$. The resulting MPC computes
\begin{equation}
    \widehat{\mathbf K}_{k,S}(\xi_k)
    \in\argmin_{\mathbf K\in\Kcal^N}
    \widehat J_{N,S}(\mathbf K;\xi_k),
    \label{eq:rhc_saa_problem}
\end{equation}
and applies the first feedback law of the optimal policy sequence
\(
    u_k=\pi_{\widehat K_{0|k,S}}(\xi_k)
\)
in a receding-horizon fashion.

After observing $x_{k+1}$, the posterior is updated using
\eqref{eq:bayes_update}.
For Bayesian linear regression, the next MPC problem uses parameter
samples from the updated Gaussian posterior:
\begin{equation}
    \theta_{k+1}^{(s)}
    =\mu_{k+1}+L_{k+1}z^{(s)},
    \qquad
    L_{k+1}L_{k+1}^\top=\Sigma_{k+1},
    \label{eq:rhc_updated_sampling_map}
\end{equation}
where $z^{(s)}\sim\mathcal N(0,I_{n_\theta})$.
The controller incorporates new observations through the updated posterior and repeated finite-horizon policy optimization.
Note that we update the posterior but do not resample the random variables \(z^{(s)}\) used to generate the model realizations.

\begin{remark}[Finite-horizon guarantee in receding horizon]
\label{rem:receding_horizon_theory}
Theorem~\ref{thm:saa_consistency} concerns the finite-horizon SAA problem at
a fixed hyperstate. In a receding-horizon implementation, this
finite-horizon problem is solved repeatedly at each hyperstate \(\xi_k\) at time $k \in \mathbb{N}_0$.
The theorem therefore justifies the sampled optimization performed at each MPC
solve, but does not ensure any closed-loop performance.
\end{remark}

\section{Numerical Results}
\label{sec:numerics}
We evaluate the proposed sampling-based dual MPC approach on two stochastic problems in Åström and Wittenmark \cite{astrom1989adaptive}. 
The computations are performed in Python using L-BFGS-B \cite{byrd1995limited} in SciPy. The code is available online at \text{\footnotesize{\url{https://github.com/TrenBaltussen/Dual-MPC-Sampling}}}.

\subsection{Example 1: Time-Invariant Unknown Gain}
Consider the stochastic integrator from \cite[Sec.~7.3]{astrom1989adaptive}
\begin{equation}
    x_{k+1}=x_k+bu_k+e_k,
\end{equation}
where $b$ is a time-invariant gain and
$e_k\sim\mathcal N(0,R_2)$ is independent Gaussian noise.
This is a scalar instance of the Bayesian linear regression
model \eqref{eq:blr_transition_model}, with
$g(x,u; b)=x+bu$.
The physical state $x_k$ and the Gaussian posterior
$\gamma_k=\mathcal N(\hat b_k,P_k)$
defines the hyperstate $\xi_k=(x_k,\hat b_k,P_k)$.
We consider the stage and terminal costs
$\ell(x,u) = \ell_N(x) = x^2$.
We use $b\sim\mathcal N(0,0.45^2)$, $R_2=0.01$, $b=0.65$, $N=8$ and $S=128$.
The Bayesian update \eqref{eq:blr_update_map}
specializes to \cite{astrom1989adaptive}
\begin{subequations}\label{eq:time_invariant_gain_update}
\begin{align}
P_{k+1}
&=\frac{P_kR_2}{R_2+P_ku_k^2},\\
\hat b_{k+1}
&=\hat b_k+
\frac{P_ku_k}{R_2+P_ku_k^2}
\left(x_{k+1}-x_k-\hat b_ku_k\right).
\end{align}
\end{subequations}
We consider an affine policy of the form
\begin{equation}
\label{eq:affine_policy}
\begin{aligned}
\psi(\xi)
&=\begin{bmatrix}x & \hat b & P & 1\end{bmatrix}^{\top},\\
K_i
&=\begin{bmatrix}a_i^1 & a_i^2 & a_i^3 & a_i^4\end{bmatrix}
\in\R^{1\times4},\\
\pi_{K_i}(\xi)&=K_i\psi(\xi).
\end{aligned}
\end{equation}
The policy parameter sequence
$\mathbf K=(K_0,\ldots,K_{N-1})\in\Kcal^N$
is shared across all $S$ sampled trajectories, cf. Section~\ref{sec:saa}.

We compare the proposed policy-based controller with four baselines that differ in how they account for model uncertainty and future observations.
\begin{itemize}
    \item Oracle: $u_k^{\mathrm{oracle}}=-\frac{1}{b_k}x_k$ has access to true model parameter and serves as a performance benchmark.
    \item Certainty equivalent: $u_k^{\mathrm{CE}} \hspace{-1mm} = \hspace{-1mm} - \frac{1}{\hat b_k}x_k$ uses the current posterior mean and has a closed-form solution \cite{astrom1989adaptive}.
    \item Open-loop SAA: $u_k=\pi_{\widehat K_{0|k,S}}(\xi_k)$ is the first input of an open-loop input sequence with $\psi(\xi) = 1$ in \eqref{eq:affine_policy}.
    \item Policy-based dual SAA: $u_k=\pi_{\widehat K_{0|k,S}}(\xi_k)$ is the first input of the optimized feedback policy of \eqref{eq:blr_saa_ocp}.
    \item Cautious: $u_k^{\mathrm{caut}}=-\frac{\hat b_k x_k}{\hat b_k^2+P_k}$ is the $N=1$-step dual controller and has a closed-form solution, cf. \cite{astrom1989adaptive}
\end{itemize}

\begin{figure}[t]
\centering
\includegraphics[width=\columnwidth]{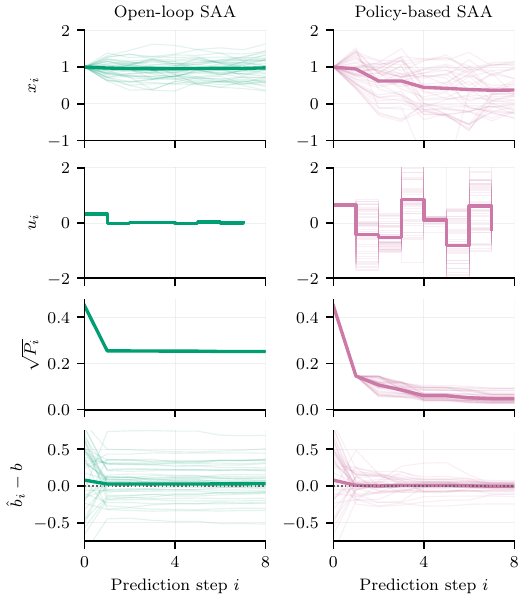}
\caption{Finite-horizon prediction for the time-invariant unknown gain problem at time $k=0$. Shaded lines show individual uncertainty realizations, solid lines show the average.}
\label{fig:finte_horizon}
\end{figure}

\begin{figure}[t]
\centering
\includegraphics[width=\columnwidth]{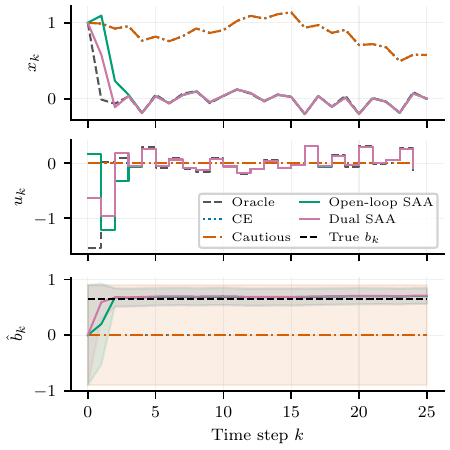}
\caption{Time-invariant gain problem. Shaded areas show the posterior standard deviation of the unknown gain $b$.}
\label{fig:fixed_gain}
\end{figure}

The coefficients \eqref{eq:affine_policy} are bounded by $|a_i^j|\leq 100$. The applied inputs are saturated at $|u_k|\leq 20$ to prevent excessive control actions, particularly for the oracle and CE controllers when $\hat b$ is close to zero. For $\hat b=0$, the $u^{\mathrm{CE}}$ is set to zero.

Fig.~\ref{fig:finte_horizon} compares the open-loop and policy-based SAA controllers on the finite-horizon problem. 
As can be seen, the different uncertainty realizations lead to different input sequences under the same closed-loop policy, while the open-loop SAA applies the same input sequence to all scenarios.
The policy-based dual controller exhibits probing behavior, resulting in a more informative posterior and steers the system to the origin more quickly than the open-loop SAA controller.
Fig.~\ref{fig:fixed_gain} shows a closed-loop experiment where the MPCs are applied in a receding-horizon fashion.
As the prior gain estimate is zero, the certainty equivalent (CE) and cautious controller apply zero input and fail to identify and regulate the system. The belief of policy-based dual MPC converges slightly faster than the open-loop SAA but both achieve similar closed-loop performance.

\subsection{Example 2: Time-Varying Unknown Gain}
To further illustrate our method, we consider the stochastic integrator with a time-varying gain from \cite[Example~7.2]{astrom1989adaptive},
\begin{equation}
\begin{aligned}
    x_{k+1} &= x_k+b_ku_k+e_k,\\
    b_{k+1} &= \varphi_b b_k+v_k,
\end{aligned}
\end{equation}
where $e_k\sim\mathcal N(0,R_2)$ and $v_k\sim\mathcal N(0,R_1)$ are independent.
Following \cite{astrom1989adaptive}, we set $\varphi_b=0.95$ and $R_1=0.09$.
This corresponds to Remark~\ref{rem:time_varying_dynamics} with $\theta_k=b_k$,
$g(x,u;b)=x+bu$, and $h(b,v)=\varphi_b b+v$.
The cost functions and policy parameterization are unchanged from Example 1.
We use $R_2=1$, $\hat b_0=0.5$, and $P_0=R_1/(1-\varphi_b^2)$.
The Gaussian belief
$\gamma_k=\mathcal N(\hat b_k,P_k)$
describes the current gain, yielding
$\xi_k=(x_k,\hat b_k,P_k)$.
The Bayesian update \eqref{eq:bayes_update} reduces to a Kalman filter for the unknown gain $b_k$ \cite[Eq.~(7.9)]{astrom1989adaptive}:
\begin{subequations}
\begin{align}
P_k^+
&=\frac{P_kR_2}{R_2+P_ku_k^2},\\
\hat b_k^+
&=\hat b_k+
\frac{P_ku_k}{R_2+P_ku_k^2}
\left(x_{k+1}-x_k-\hat b_ku_k\right),\\
\hat b_{k+1}
&=\varphi_b\hat b_k^+,\\
P_{k+1}
&=\varphi_b^2P_k^++R_1.
\end{align}
\end{subequations}
These equations provide an exact finite-dimensional
Bayesian update for the sampled hyperstate trajectories.

\begin{figure}[t]
\centering
\includegraphics[width=\columnwidth]{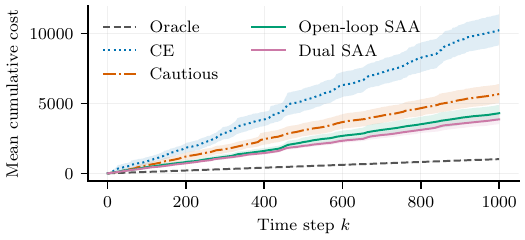}
\caption{Closed-loop cost in the time-varying gain problem over 20 Monte Carlo runs with 95\% confidence intervals.}
\label{fig:cost_varying_gain}
\end{figure}

\begin{table}[t]
\centering
\caption{Mean and standard deviation (SD) of 20 Monte Carlo runs for the time-varying gain example for $T=1000$.}
\label{tab:varying_gain_performance}
\small
\setlength{\tabcolsep}{4pt}
\begin{tabular}{@{}l*{6}{r}@{}}
\toprule
& \multicolumn{2}{c}{Normalized regret} &
\multicolumn{2}{c}{Gain RMSE} &
\multicolumn{2}{c}{Input energy} \\
Controller & Mean & SD & Mean & SD & Mean & SD \\
\midrule
Oracle        & 0.00 & 0.00 & \multicolumn{2}{c}{---} & 26.3 & 4.8 \\
Dual SAA      & 1.00 & 0.32 & 0.567 & 0.040 & 3.44 & 1.04 \\
OL SAA        & 1.16 & 0.44 & 0.597 & 0.039 & 2.56 & 0.85 \\
Cautious      & 1.64 & 0.54 & 0.641 & 0.046 & 2.43 & 1.02 \\
CE            & 3.24 & 0.91 & 0.475 & 0.033 & 58.6 & 11.6 \\
\bottomrule
\end{tabular}
\end{table}

\begin{figure*}[t]
\centering
\begin{subfigure}[t]{0.49\textwidth}
    \centering
    \includegraphics[width=\linewidth]{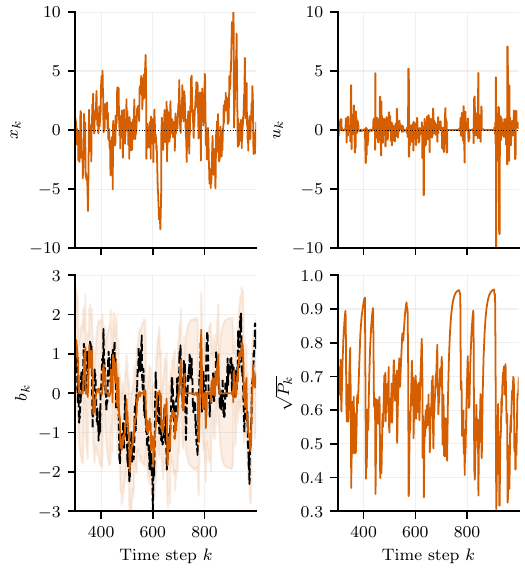}
    \caption{Cautious controller.}
    \label{fig:cautious_varying_gain}
\end{subfigure}\hfill
\begin{subfigure}[t]{0.49\textwidth}
    \centering
    \includegraphics[width=\linewidth]{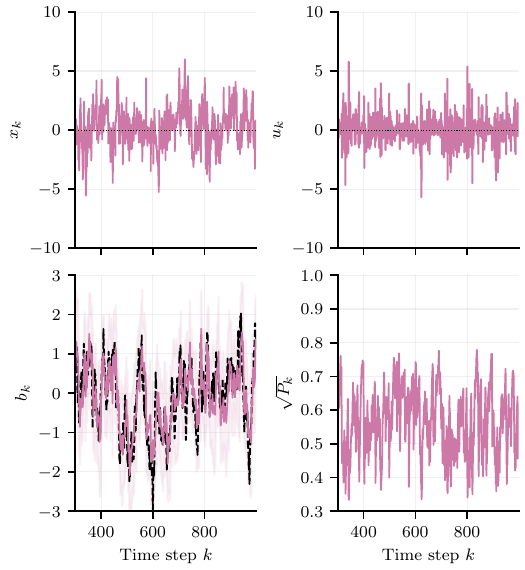}
    \caption{Policy-based SAA Dual MPC.}
    \label{fig:dual_mpc_varying_gain}
\end{subfigure}
\caption{Representative closed-loop simulations in the time-varying gain example for $k=300$--$1000$. Black dashed lines indicate true parameter $b_k$, colored line indicates posterior mean $\hat b_k$, and shaded area indicates posterior standard deviation.}
\label{fig:varying_gain_zoom_comparison}
\end{figure*}

The results of the closed-loop Monte Carlo simulations are summarized in Fig.~\ref{fig:cost_varying_gain} and Table~\ref{tab:varying_gain_performance}. 
Regret is the difference between the cost of the controller and the cost of the oracle controller, the root mean square error (RMSE) of the gain is computed between the true and estimated gains, and the input energy $(u^2)$ is averaged over time.
Fig.~\ref{fig:cautious_varying_gain} and Fig.~\ref{fig:dual_mpc_varying_gain} show representative closed-loop simulations of the cautious and proposed dual SAA controllers, respectively.

The CE controller shows the highest regret and input energy as a result of poor parameter estimation and aggressive control actions.
The cautious controller exhibits the so-called \emph{turn-off behavior} \cite{astrom1989adaptive} around time steps 400, 750, and 850, where the estimated gain approaches zero while the posterior variance increases, cf. Fig.~\ref{fig:cautious_varying_gain}. This leads to a vicious cycle of substantial reduction in input magnitude and limited information obtained from subsequent measurements.
The cautious controller accounts for parameter uncertainty, leading to lower input energy than the CE controller, but it exhibits higher regret and gain RMSE than the proposed dual SAA, potentially, due to its myopic planning and turn-off behavior.
The open-loop SAA exhibits lower regret and gain RMSE than the cautious controller, but it does not actively probe the system to acquire information about the unknown gain.
Finally, the proposed dual SAA achieves the lowest regret and gain RMSE among the controllers. As the dual SAA controller can select probing inputs, the input energy is higher than the cautious and open-loop SAA controllers, but significantly lower than the CE controller, indicating that the dual MPC carefully selects probing inputs.
The proposed dual MPC achieves the lowest cumulative tracking cost among the controllers, demonstrating its effectiveness in intrinsically balancing exploration and exploitation in the presence of stochastic, time-varying uncertainty.

\section{Conclusion}

This paper proposes a sampling-based MPC approach for dual control. The method approximates the exact finite-horizon Bayesian objective using a sample average approximation and optimizes over hyperstate-feedback policies. The resulting controller propagates sample-dependent future beliefs while keeping the decision variables tied across samples through a shared policy sequence. 
The proposed method is a practical approach to dual control, and the numerical examples illustrate its ability to intrinsically balance exploration and exploitation in the presence of model uncertainty.
Extensions toward nonparametric models were discussed, and are considered part of future work. Investigating closed-loop properties, such as optimality, stability and constraint satisfaction, is an important direction for future research.

\bibliographystyle{IEEEtran}
\bibliography{references}

\end{document}